\documentclass[12pt]{article}
\usepackage{amsmath}
\usepackage{amssymb}
\usepackage{amsopn}
\usepackage{latexsym}
\usepackage{amsfonts}
\title{\textbf{The Armstrong-Frederick cyclic hardening plasticity model with {C}osserat effects}}
\author{\textbf{Krzysztof Che{\l}mi\'{n}ski}\\[0.5ex]
\textbf{\footnotesize{Faculty of Mathematics and
 Information Science, Warsaw University of Technology,}}\\[-1ex]
\textbf{\footnotesize{pl. Politechniki 1, 00-661 Warsaw, Poland}}\\[-1ex]
\textbf{\footnotesize{E-Mail: kchelmin@mini.pw.edu.pl}}\\[1ex]
\textbf{Patrizio Neff}\\[0.5ex]
\textbf{\footnotesize{Fakult\"at f\"ur Mathematik, Universit\"at
    Duisburg-Essen,}}\\[-1ex] 
\textbf{\footnotesize{Lehrstuhl f\"ur Nichtlineare Analysis und Modellierung}}\\[-1ex]
\textbf{\footnotesize{Thea-Leymann Strasse 9,  45141 Essen, Germany}}\\[-1ex]
\textbf{\footnotesize{E-Mail: patrizio.neff@uni-due.de}}\\[1ex]
\textbf{Sebastian Owczarek}\\[0.5ex] 
\textbf{\footnotesize{Faculty of Mathematics and
 Information Science, Warsaw University of Technology,}}\\[-1ex]
\textbf{\footnotesize{pl. Politechniki 1, 00-661 Warsaw, Poland}}\\[-1ex]
\textbf{\footnotesize{E-Mail: s.owczarek@mini.pw.edu.pl}}}
\date{}
\newtheorem{tw}{Theorem}[section]
\newtheorem{lem}[tw]{Lemma}
\newtheorem{de}[tw]{Definition}
\DeclareMathOperator{\dev}{dev}
\DeclareMathOperator{\so}{\mathfrak{so}}
\begin{document}
\maketitle
\begin{abstract}
\noindent
We propose an extension of the cyclic hardening plasticity model formulated by Armstrong and Frederick which includes micropolar effects. Our micropolar extension establishes coercivity of the model which is otherwise not present. We study then existence of solutions to the quasistatic, rate-independent Armstrong-Frederick model with Cosserat effects which is, however, still of non-monotone, non-associated type. In order to do this, we need to relax the pointwise definition of the flow rule into a suitable weak energy-type inequality.  It is shown that the limit in the Yosida approximation process satisfies this new solution concept. The limit functions have a better regularity than previously known in the literature, where the original Armstrong-Frederick model has been studied.
\end{abstract}
\newcommand{\nn}{\nonumber}
\newcommand{\KK}{\sigma_{\rm y}}
\newcommand{\ve}{\varepsilon}
\newcommand{\R}{{\mathbb R}}
\newcommand{\D}{{\mathbb C}}
\newcommand{\E}{{\cal E}}
\newcommand{\K}{{\cal K}}
\renewcommand{\S}{{\cal S}^3}
\renewcommand{\SS}{{\cal S}^3_{\dev}}
\newcommand{\id}{ {1\!\!\!\:1 } }
\section{Introduction}
One of the well-known models to describe cyclic plasticity is the non-linear kinematic hardening model formulated by Armstrong and Frederick \cite{3}. This model has been highly rated, because it is based on a physical mechanism of strain hardening and dynamic recovery, and because it has the capability of representing reasonably well the shapes of stress-strain hysteresis loops, especially those of constant strain ranges. Therefore, implementation of the Armstrong-Frederick model in finite element methods has been examined in several studies to date. Thus, that model is now available as a material model of cyclic hardening plasticity in commercial, general-purpose software for finite element analysis.

The Armstrong-Frederick model (AF) is a modification of the Melan-Prager model, which is well know in the literature and it can also be seen as an approximation of the Prandtl-Reuss model. The key modification of this simple model is the so-called "recall"-term, changing the evolution law for the symmetric backstress tensor $b$ from a classical linear kinematic hardening law (Melan-Prager) to a nonlinear kinematic hardening law, i.e.,
\begin{align*}
  b_t=\underbrace{c\, \varepsilon_t^p}_{\text{lin. kin. hardening}}-\underbrace{d\, |\varepsilon_t^p| b}_{\text{recall-term, nonlinear hardening}}\, .
\end{align*}
Here, $\varepsilon^p$ is the symmetric plastic strain tensor, $c$ and $d$ are positive material constants. The "recall"-term entails the $L^{\infty}$-boun\-ded\-ness of the backstress $b$, a property which is an experimental fact since to the contrary, for high frequency cycles softening and rupture will occur. Therefore, the AF-model shows nonlinear kinematic hardening, but only to within a certain extent. The more realistic description of cyclic hardening plasticity experiments with the AF-model, however, has a prize to be paid: the model is non-coercive (bounded hardening), it is of non-monotone type and not of gradient type (non-associated flow rule). Thus, the AF-model is one of the prominent small strain plasticity models which has yet escaped the efforts of mathematicians to establish well-posedness. 

The mathematical analysis being quite challenging, there are no encompassing existence results for this model in the literature. The first (partial) mathematical result for the Armstrong-Frederick model was obtained by the first author in the article \cite{7}. There, the non-monotone, non-associated AF-model was written as a model of pre-monotone type (for the definition we refer to \cite{2}). In this work the author used a Yosida approximation to the monotone part of the flow rule. The obtained a priori estimates are, however, not sufficient to pass to the limit with such approximations and to get $L^2$-strong solutions (see Section 3 in \cite{7}). It was only shown that the limit functions satisfy the so called "reduced energy inequality". In the article \cite{14} a regularization of the "recall"-term in the equation for the backstress was proposed. The existence of a rescaled in time solution to the Armstrong-Frederick model with the regularized equation for the backstress could then be established. The rescaling idea is very simple: a new time variable $\tau=\zeta(t)$ is proposed. Then the new system is easier to analyze, because the plastic strain is now uniformly Lipschitz with respect to the rescaled time. The main problem is to get back to the original system with the rescaled in time solution. It is, in principle, possible for rate-independent models, where the flow rule is invariant under scaling of the time variable. The Armstrong-Frederick model is rate-independent but the authors of \cite{14} are not able to get back to the original system. The rescaling idea has already been applied in the plasticity context in \cite{4,12,13}.

In this paper we want to extend the system of equations proposed by Armstrong and Frederick to include micropolar effects. In the classical metal perfect plasticity models at infinitesimal strains it has been shown in a series of papers \cite{8,9,23,24} that a coupling with Cosserat elasticity may also regularize the ill-posedness of the Prandtl-Reuss plasticity model. This is possible because the Cosserat coupling leads to coercivity. Perfect plasticity, however, is yet characterized by a monotone flow rule of gradient type (associated plasticity). Therefore, the question arises naturally, whether adding Cosserat microrotations to the AF-model is still enough to regularize the problem in the way to satisfy the flow rule in a standard pointwise sense. From a modelling perspective, adding
microrotations means to consider a material made up of individual particles which can rotate and interact with each other \cite{16,17,22,18}. For phenomenological polycrystalline plasticity adding Cosserat effects is arguably a physically motivated regularization: the individual crystal grains are rotating and interacting with each other. 

The extension of the Armstrong-Frederick model to include Cosserat effects follows the lines proposed in \cite{8}, where the authors added the Cosserat effect to the classical elastoplasticity model with a monotone flow rule. There, and in our present approach, only the elasticity relation is augmented with Cosserat effects, the plastic constitutive equations, and notably the "recall"-term, is left unchanged, contrary to \cite{14}. Regarding the effect of the Cosserat-modification for classical plasticity models, it has been proved that the new model is thermodynamically admissible and that there exists a unique, global in time solution to Cosserat elasto-plasticity. In \cite{10} a $H^1_{loc}$- regularity result for the stresses and strains was proved, cf. \cite{19}. The dynamic Cosserat plasticity was studied in \cite{9}, see also \cite{23,24}.
Another way to regularize classical plasticity models is by introducing gradient plasticity effects \cite{29,25,26,27,28}. 
However, a modification of the AF-model to include higher gradients will be left to future work. Moreover, the coupling with thermal effects can be treated as another attempt to regularize models from the inelastic deformation theory (c.f. \cite{31,32}).

Many non-monotone models from the theory of inelastic deformation processes in metals are also non-coercive (for the definition see \cite{2}) and the existence results for such models is only very weakly examined. For example: the solutions obtained in articles \cite{7}, and \cite{14} had a low regularity with respect to time and space (see also \cite{20, 30}, where the non-monotone model of poroplastisity was considered). In our opinion, it is expedient to 
first consider non-monotone but coercive models and to obtain better regularity results for the solutions. This article presents the first mathematical result in this respect for our new Armstrong-Frederick model with Cosserat effects, which is non-monotone, non-associated but coercive. 


\section{The Armstrong-Frederick model with Cosserat effect}
This section is devoted to the formulation of the Armstrong-Frederick model with Cosserat effects.

From the mechanical results for Cosserat plasticity (see for example \cite{8}, \cite{9}) we conclude that we deal with the following initial-boundary value problem: we are looking for the displacement field $u:\Omega\times[0,T]\rightarrow \R^3$, the microrotation matrix
$A:\Omega\times[0,T]\rightarrow\mathfrak{so}(3)$ ($\mathfrak{so}(3)$ is the set of skew-symmetric $3\times 3$ matrices) and the vector of internal variables\\ $z=(\ve^p,b):\Omega\times [0,T]\rightarrow\SS\times\SS$ ($\ve^p$ is the classical infinitesimal symmetric plastic strain tensor, $b$ is the symmetric backstress tensor and the space $\SS$ denotes the set of symmetric $3\times 3$-matrices with vanishing trace) satisfying the following system of equations
\renewcommand{\theequation}{\thesection.\arabic{equation}}
\setcounter{equation}{0}%
\begin{eqnarray}
\label{eq:2.1}
\mathrm{div}_x T&=&-f\,,\nn\\
T&=&2\mu(\ve(u)-\ve^p)+2\mu_c(\mathrm{skew}(\nabla_x u)-A)+\lambda\mathrm{tr}(\ve(u)-\ve^p)\id\,,\nn\\[1ex]
-l_c\,\Delta_x\mathrm{axl}\,(A)&=&\mu_c\,\mathrm{axl}\,(\mathrm{skew}(\nabla_x u)-A)\,,\nn\\[1ex]
\ve^{p}_{t}&\in&\partial I_{K(b)}\Big (T_{E}\Big)\,,\\
T_E&=&2\mu(\ve(u)-\ve^p)+\lambda\mathrm{tr}(\ve(u)-\ve^p)\id,\nn\\
b_t&=&c\,\ve^p_t-d\,|\ve_t^p|b\,,\nn
\end{eqnarray}
where $\ve(u)=\mathrm{sym}(\nabla_x u)$ denotes the symmetric part of the gradient of the displacement. The above equations are studied for $x\in\Omega\subset \R^3$ and $t\in [0,T]$, where $\Omega\subset\R^3$ is a bounded domain with smooth
boundary $\partial\Omega$ and $t$ denotes the time.\\
The set of admissible elastic stresses $K(b(x,t))$ is defined in the form\\
$K(b)=\{T_E\in\S\,:\,|\dev\,(T_E)-b|\leq \KK\}$, where $\dev\,(T_E)=T_E-\frac{1}{3}\mathrm{tr}(T_E)\cdot\id$, $\KK$ is a material parameter (the yield limit) and $\id$ denotes the identity matrix. The function $I_{K(b)}$ is the indicator function of the admissible set $K(b)$ and $\partial I_{K(b)}$ is the subgradient of the convex, proper, lower semicontinous function $I_{K(b)}$. The function $f:\Omega\times[0,T]\rightarrow \R^3$ describes the density of the applied body forces, the parameters $\mu$, $\lambda$ are positive Lam\'e constants (the elastic constitutive equation
can be generalized in the obvious way to anisotropic case), $\mu_c>0$ is the Cosserat couple modulus and $l_c>0$ is a material parameter with dimensions $[m^2]$, describing a length scale of the model due to the Cosserat effects.  $c,d>0$ are material constants. The operator $\mathrm{skew}(T)=\frac{1}{2}(T-T^T)$ denotes the skew-symmetric part of a $3\times 3$-tensor. The operator $\mathrm{axl}:\mathfrak{so}(3)\rightarrow \R^3$ establishes the identification of a skew-symmetric matrix with vectors in $\R^3$. This means that if we take $A\in\mathfrak{so}(3)$, which is in the form $A=((0,\alpha,\beta),(-\alpha,0,\gamma),(-\beta,-\gamma,0))$, then $\mathrm{axl}(A)=(\alpha,\beta,\gamma)$.\\
Notice again that the expression $|\ve^p_t|b$ is a perturbation of Melan-Prager model - if $d=0$ then we obtain the Melan-Prager linear kinematic hardening model.\\ 
The system (\ref{eq:2.1}) is considered with mixed boundary conditions for the displacement:
\begin{eqnarray}
\label{eq:2.2}
u(x,t)&=&g_D(x,t)\qquad \textrm{ for}\quad x\in\Gamma_D \quad\textrm{and}\quad t\geq 0\textrm{,}\nn\\
T(x,t)\cdot n(x)&=&g_N(x,t)\qquad\; \textrm{for}\quad x\in\Gamma_N \quad\textrm{and}\quad t\geq 0\textrm{,}
\end{eqnarray}
where $n(x)$ is the exterior unit normal vector to the boundary $\partial\Omega$ at the point $x$, $\Gamma_D$, $\Gamma_N$ are open subsets of $\partial\Omega$ satisfying $\partial\Omega= \bar{\Gamma}_D \cup \bar{\Gamma}_N$, $\Gamma_D\cap\Gamma_N=\emptyset$ and $\mathcal{H}_2(\Gamma_D)>0$- $\mathcal{H}_2$ denotes the $2$-dimensional Hausdorff measure. The system (\ref{eq:2.1}) is considered with Dirichlet boundary conditions for the microrotation:
\begin{eqnarray}
\label{eq:2.3}
A(x,t)&=&A_D(x,t)\; \textrm{ for}\quad x\in \partial \Omega \quad\textrm{and}\quad  t\geq 0.
\end{eqnarray}
Finally, we consider the system (\ref{eq:2.1}) with the following initial conditions
\begin{eqnarray}
\label{eq:2.4}
\ve^p(x,0)=\ve^{p,0}(x),\qquad b(x,0)=b^0(x).
\end{eqnarray}
The free energy function associated with the system (\ref{eq:2.1}) is given by the formula
\begin{eqnarray}
\label{eq:2.5}
\rho\psi(\ve,\ve^p,A,b)&=&\mu\|\ve(u)-\ve^p\|^2+
\mu_c\|\mathrm{skew}(\nabla_x u)-A\|^2\nn\\
&+&\frac{\lambda}{2}\Big(\mathrm{tr}(\ve(u)-\ve^p)\Big)^2+2l_c\|\nabla_x\mathrm{axl}(A)\|^2+\frac{1}{2c}\|b\|^2,
\end{eqnarray}
where $\rho$ is the mass density which we assume to be constant in time and space. The total energy is of the form:
\begin{eqnarray*}
\E(\ve,\ve^p,A,b)(t)=\int\nolimits_{\Omega}\rho\psi(\ve(x,t),\ve^p(x,t),A(x,t),b(x,t))dx.
\end{eqnarray*}
From Section $2$ of the article \cite{7} we know that the inelastic constitutive equation occuring in (\ref{eq:2.1}) is of pre-monotone type. Moreover, if initial data $(\ve^{p,0},b^0)\in\SS\times\SS$, then any solution $(\ve^p(t),b(t))$ (if there exists) belongs to $\SS\times\SS$, because the right hand side of $(\ref{eq:2.1})_4$ is a subset of $\SS$ and the following system
\begin{eqnarray*}
\frac{d}{dt}(\mathrm{tr}\,\ve^p) &=&0 \quad \mathrm{with}\quad \mathrm{tr}\,\ve^p(0)=0\, ,\\ 
\frac{d}{dt}(\mathrm{tr}\,b)+d|\ve^p_t|\mathrm{tr}\, b &=&0 \quad \mathrm{with}\quad \mathrm{tr}\,b(0)=0\, ,
\end{eqnarray*}
possesses the unique solution $(\mathrm{tr}\,\ve^p,\mathrm{tr}\, b)=(0,0)$. Inspired by the work \cite{7}, we propose to rewrite the set of admissible stresses in the form
$$\K=\{(\dev\,(T_E),b)\in\SS\times\SS\,:\, |\dev\,(T_E)+cb|\leq \KK\},$$
where the constant $\KK$ is the same as in the yield condition and the inelastic constitutive equation can now be written in the equivalent form
\begin{displaymath}
\frac{d}{dt}\left(\begin{array}{c}
\ve^p\\
b\end{array}\right) \in \partial I_{\K}\Big((\dev\,(T_{E}),-\frac{1}{2}b)\Big)-\left(\begin{array}{c}
0\\
d\,|\ve^p_t|b\end{array}\right).
\end{displaymath}
\section{Main result}
Here we define a notion of the definition of the solution for the system (\ref{eq:2.1}). Next we formulate the main result of this paper. First, let us start with the definition of $L^2$-strong solution of the system (\ref{eq:2.1}).
\renewcommand{\theequation}{\thesection.\arabic{equation}}
\setcounter{equation}{0}%
\begin{de}
\label{de:3.1}
Fix $T>0$. We say that a vector $(u,A,T,\ve^p,b)\in W^{1,\infty}(0,T;H^1(\Omega;\R^3)\times H^2(\Omega;\mathfrak{so}(3))\times L^2(\Omega;\R^9)\times L^2(\Omega;\SS)\times L^2(\Omega;\SS))$ is $L^2$-strong solution of the system
\begin{eqnarray}
\mathrm{div}_x T&=&-f\,,\nn\\
T&=&2\mu(\ve(u)-\ve^p)+2\mu_c(\mathrm{skew}(\nabla_x u)-A)+\lambda\mathrm{tr}(\ve(u)-\ve^p)\id\,,\nn\\[1ex]
-l_c\,\Delta_x\mathrm{axl}\,(A)&=&\mu_c\,\mathrm{axl}\,(\mathrm{skew}(\nabla_x u)-A)\,,\nn\\[2ex]
\frac{d}{dt}\left(\begin{array}{c}
\ve^p\\
b\end{array}\right)+\left(\begin{array}{c}
0\\
d\,|\ve^p_t|b\end{array}\right) &\in& \partial I_{\K}\left(\begin{array}{c}
\dev\,(2\mu(\ve(u)-\ve^p))\\
-\frac{1}{c}\,b\end{array}\right) \nn
\end{eqnarray}
if
\begin{enumerate}
\item $|\ve^p_t|b\in L^{\infty}(0,T;L^2(\Omega,\SS))$,
\item $|\dev\,\Big(2\mu(\ve(u(x,t))-\ve^p(x,t))\Big)-b(x,t)|\leq \KK$ for almost all $(x,t)\in \Omega\times (0,T)$,
\item the equations and inclusion above are satisfied for almost all $(x,t)\in \Omega\times (0,T)$.
\end{enumerate}
\end{de}
Next, we are going to define a weaker notion of solutions to the system (\ref{eq:2.1}). We give a motivation of this definition. It will be the energy inequality combined with special test functions.\\
Let us consider another convex set (which will be used as set of test functions further on)
$$\K^{\ast}=\{(\dev\,(T_E),-\frac{1}{c}b)\in\SS\times\SS\,:\, |\dev\,(T_E)-b|+\frac{d}{2c}|b|^2\leq \KK\},$$
where the constant $\KK$ is the same as in the yield condition. It is not difficult to see that for all $L^2$-strong solutions we get (for details we refer to figure 1 and 2 of \cite{7})
\begin{equation}
\label{eq:3.1}
\ve^p_t\Big(\dev\,(T_E)-\dev\,(\hat{T}_E)\Big)+b_t\Big(-\frac{1}{c}\,b+\frac{1}{c}\,\hat{b}\Big)\geq 0\;\;\;\mathrm{for\; all}\;\; (\dev\,(\hat{T}_E),\hat{b})\in \K^{\ast}\,,\quad
\end{equation}
where $\hat{T}_E\in L^{\infty}(0,T;L^2(\Omega,\S))$ is any stress tensor. From the theory of elasticity, we know that there exists a positive definite operator $\D^{-1}:\S\rightarrow\S$ such that $\D^{-1}T_{E,t}=\ve_t-\ve^p_t$. Integrating (\ref{eq:3.1}) over $\Omega\times (0,t)$ for $t\leq T$ we obtain 
\begin{eqnarray}
\label{eq:3.2}
&&\int\nolimits_0^t\int\nolimits_{\Omega}\D^{-1}T_{E,t}(x,\tau)T_E(x,\tau)dxd\tau+ \frac{1}{c}\int\nolimits_0^t\int\nolimits_{\Omega}b_t(x,\tau)b(x,\tau)dxd\tau\nn\\[1ex]
&\leq&\int\nolimits_0^t\int\nolimits_{\Omega}\ve_t(x,\tau)(T_E(x,\tau)-\hat{T}_E(x,\tau))dxd\tau+ \int\nolimits_0^t\int\nolimits_{\Omega}\D^{-1}T_{E,t}(x,\tau)\hat{T}_E(x,\tau)dxd\tau\nn\\[1ex]
&+& \frac{1}{c}\int\nolimits_0^t\int\nolimits_{\Omega}b_t(x,\tau)\hat{b}(x,\tau)dxd\tau= \int\nolimits_0^t\int\nolimits_{\Omega}\nabla u_t(x,\tau)(T(x,\tau)-\hat{T}_E(x,\tau))dxd\tau\nn\\[1ex]
&-&2\mu_c\int\nolimits_0^t\int\nolimits_{\Omega}\Big(\mathrm{skew}(\nabla_x u(x,\tau))-A(x,\tau)\Big)\mathrm{skew}(\nabla_x u_t(x,\tau))dxd\tau\nn\\[1ex]
&+& \frac{1}{c}\int\nolimits_0^t\int\nolimits_{\Omega}b_t(x,\tau)\hat{b}(x,\tau)dxd\tau\,.
\end{eqnarray}
Let us assume that $\mathrm{div}\; \hat{T}_E\in L^{\infty}(0,T;L^2(\Omega,\R^3))$ (now the trace of $\hat{T}_E$ in the normal direction to the boundary $\partial\Omega$ is well defined) and $\hat{T}_E(x,t)\cdot n(x)=g_N(x,t)$ on $\Gamma_N\subset\partial\Omega$. Integrating by parts in the first term on the right hand side of (\ref{eq:3.2}), using equation $(\ref{eq:2.1})_1$ and the boundary data we have the following inequality
\begin{eqnarray}
\label{eq:3.3}
&&\frac{1}{2}\int\nolimits_{\Omega}\D^{-1}T_{E}(x,t)T_E(x,t)dx+ \mu_c\int\nolimits_{\Omega}|\mathrm{skew}(\nabla_x u(x,t))-A(x,t)|^2dx+ \frac{1}{2c}\int\nolimits_{\Omega}|b(x,t)|^2dx\nn\\[1ex]
&\leq&\frac{1}{2}\int\nolimits_{\Omega}\D^{-1}T_{E}(x,0)T_E(x,0)dx+ \mu_c\int\nolimits_{\Omega}|\mathrm{skew}(\nabla_x u(x,0))-A(x,0)|^2dx+ \frac{1}{2c}\int\nolimits_{\Omega}|b(x,0)|^2dx\nn\\[1ex] &+&\int\nolimits_0^t\int\nolimits_{\Omega}u_t(x,\tau)f(x,\tau)dxd\tau+ \int\nolimits_0^t\int\nolimits_{\Omega}u_t(x,\tau)\mathrm{div}\,\hat{T}_E(x,\tau))dxd\tau\\[1ex]
&+&\int\nolimits_0^t\int\nolimits_{\Gamma_D}g_{D,t}(x,\tau)(T(x,\tau)-\hat{T}_E(x,\tau))\cdot n(x)dSd\tau+ \int\nolimits_0^t\int\nolimits_{\Omega}\D^{-1}T_{E,t}(x,\tau)\hat{T}_E(x,\tau)dxd\tau\nn\\[1ex]
&+& \frac{1}{c}\int\nolimits_0^t\int\nolimits_{\Omega}b_t(x,\tau)\hat{b}(x,\tau)dxd\tau
-2\mu_c\int\nolimits_0^t\int\nolimits_{\Omega}\Big(\mathrm{skew}(\nabla_x u(x,\tau))-A(x,\tau)\Big)A_t(x,\tau)dxd\tau.\nn
\end{eqnarray}
Using $(\ref{eq:2.1})_3$ to the last term of the above inequality and integrating by parts in the last term on the right hand side of (\ref{eq:3.3}) we obtain
\begin{eqnarray}
\label{eq:3.4}
&&\frac{1}{2}\int\nolimits_{\Omega}\D^{-1}T_{E}(x,t)T_E(x,t)dx+ \mu_c\int\nolimits_{\Omega}|\mathrm{skew}(\nabla_x u(x,t))-A(x,t)|^2dx\nn\\[1ex]
&+&2l_c\int\nolimits_{\Omega}|\nabla\mathrm{axl}(A(x,t))|^2dx+ \frac{1}{2c}\int\nolimits_{\Omega}|b(x,t)|^2dx
\leq\frac{1}{2}\int\nolimits_{\Omega}\D^{-1}T_{E}(x,0)T_E(x,0)dx\nn\\[1ex]
&+& \mu_c\int\nolimits_{\Omega}|\mathrm{skew}(\nabla_x u(x,0))-A(x,0)|^2dx+ \frac{1}{2c}\int\nolimits_{\Omega}|b(x,0)|^2dx + 2l_c\int\nolimits_{\Omega}|\nabla\mathrm{axl}(A(x,0))|^2dx\nn\\[1ex]
 &+&\int\nolimits_0^t\int\nolimits_{\Omega}u_t(x,\tau)f(x,\tau)dxd\tau+ \int\nolimits_0^t\int\nolimits_{\Omega}u_t(x,\tau)\mathrm{div}\;\hat{T}_E(x,\tau))dxd\tau\\[1ex]
&+&\int\nolimits_0^t\int\nolimits_{\Gamma_D}g_{D,t}(x,\tau)(T(x,\tau)-\hat{T}_E(x,\tau))\cdot n(x)dSd\tau+ \int\nolimits_0^t\int\nolimits_{\Omega}\D^{-1}T_{E,t}(x,\tau)\hat{T}_E(x,\tau)dxd\tau\nn\\[1ex]
&+& \frac{1}{c}\int\nolimits_0^t\int\nolimits_{\Omega}b_t(x,\tau)\hat{b}(x,\tau)dxd\tau 
+4l_c\int\nolimits_0^t\int\nolimits_{\partial\Omega}\nabla\mathrm{axl}(A(x,\tau))\cdot n \;\mathrm{axl}(A_{D,t}(x,\tau))dSd\tau,\nn
\end{eqnarray}
where the boundary integrals are defined in the sense of duality between the space 
$$H^{\frac{1}{2}}(\partial\Omega;\R^3)\;\;\;\mathrm{ and }\;\;\; H^{-\frac{1}{2}}(\partial\Omega;\R^3)\quad \mathrm{(see\; \cite{1}\; for\; details)}.$$
$\textrm{ }$\\
Let us  assume that for all $T>0$ the given data $F$, $g_D$, $g_N$, $A_D$ have the regularity
\begin{equation}
\label{eq:3.5}
f\in H^{1}(0,T;L^2(\Omega;\R^3)),\quad g_D\in H^{1}(0,T;H^{\frac{1}{2}}(\Gamma_D;\R^3))\textrm{,}
\end{equation}
\begin{equation}
\label{eq:3.6}
g_N\in H^{1}(0,T;H^{-\frac{1}{2}}(\Gamma_N;\R^3)), \quad A_D\in H^{1}(0,T;H^{\frac{3}{2}}(\partial\Omega;\so(3))).
\end{equation}
Additionally let us assume that the initial data $(\ve^{p,0},b^0)\in L^2(\Omega;\SS)\times L^2(\Omega;\SS)$\\ satisfy
\begin{equation}
\label{eq:3.7}
|b^{0}(x)|\leq\frac{c}{d}\quad \mathrm{and}\quad |\dev\,(T^0_E(x))-b^0(x)|\leq \KK \quad \mathrm{for\; almost\; all\;}x\in\Omega,
\end{equation}
where the initial stress $T^0_E=2\mu(\ve(u(0))-\ve^{p,0})+\lambda\mathrm{tr}(\ve(u(0))-\ve^{p,0})\id\in L^2(\Omega;\S)$ is the unique solution of the following linear problem
\begin{eqnarray}
\label{eq:3.8}
\mathrm{div}_x T^0(x)&=&-f(x,0),\nn\\[1ex]
-l_c\,\Delta_x\mathrm{axl}\,(A(x,0))&=&\mu_c\,\mathrm{axl}\,(\mathrm{skew}(\nabla_x u(x,0))-A(x,0)),\nn\\[1ex]
u(x,0)_{|_{\Gamma_D}}=g_D(x,0) && T^0(x)\cdot n_{|_{\Gamma_N}}=g_N(x,0)\quad A(x,0)_{|_{\partial\Omega}}=A_D(x,0),
\end{eqnarray}
with
$$T^0(x)=2\mu\,(\ve(u(x,0))-\ve^{p,0}(x))+2\mu_c\,(\mathrm{skew}(\nabla_x u(x,0))-A(x,0))+\lambda\,\mathrm{tr}(\ve(u(x,0))-\ve^{p,0}(x))\id.\nn\\$$
\begin{de}
\label{de:3.2}
$\mathrm{(solution\; concept-energy\; inequality)}$\\
Fix $T>0$. Suppose that the given data satisfy (\ref{eq:3.5}) and (\ref{eq:3.6}). We say that a vector $(u,T,A,\ve^p,b)\in L^{\infty}(0,T;H^{1}(\Omega;\R^3)\times L^2(\Omega;\S)\times H^2(\Omega;\mathfrak{so}(3))\times (L^{\infty}(\Omega;\SS))^2)$ solves the problem (\ref{eq:2.1})-(\ref{eq:2.4}) if
$$(u_t,T_t,A_t,\ve^p_t,b_t)\in L^{2}(0,T;H^{1}(\Omega;\R^3)\times L^2(\Omega;\R^9)\times H^2(\Omega;\mathfrak{so}(3))\times (L^{2}(\Omega;\SS))^2),$$
the equations $(\ref{eq:2.1})_1$ and $(\ref{eq:2.1})_3$ are satisfied pointwise almost everywhere on $\Omega\times (0,T)$ and for all test functions $(\hat{T}_E,\hat{b})\in L^{2}(0,T;L^{2}(\Omega;\S)\times L^{2}(\Omega;\SS))$ such that 
$$ (\dev\,(\hat{T}_E),\hat{b})\in \K^{\ast},\quad \mathrm{div}\; \hat{T}_E\in L^{2}(0,T;L^2(\Omega,\R^3)),\quad \hat{T}_E\cdot n=g_N\; \mathrm{on}\; \Gamma_N\times (0,T),$$
the inequality
\begin{eqnarray}
\label{eq:3.9}
&&\frac{1}{2}\int\nolimits_{\Omega}\D^{-1}T_{E}(x,t)T_E(x,t)dx+ \mu_c\int\nolimits_{\Omega}|\mathrm{skew}(\nabla_x u(x,t))-A(x,t)|^2dx\nn\\[1ex]
&+&2l_c\int\nolimits_{\Omega}|\nabla\mathrm{axl}(A(x,t))|^2dx+ \frac{1}{2c}\int\nolimits_{\Omega}|b(x,t)|^2dx
\leq\frac{1}{2}\int\nolimits_{\Omega}\D^{-1}T^0_{E}(x)T^0_E(x)dx\nn\\[1ex]
&+& \mu_c\int\nolimits_{\Omega}|\mathrm{skew}(\nabla_x u(x,0))-A(x,0)|^2dx+ \frac{1}{2c}\int\nolimits_{\Omega}|b(x,0)|^2dx+2l_c\int\nolimits_{\Omega}|\nabla\mathrm{axl}(A(x,0))|^2dx\nn\\[1ex]
 &+&\int\nolimits_0^t\int\nolimits_{\Omega}u_t(x,\tau)f(x,\tau)dxd\tau+ \int\nolimits_0^t\int\nolimits_{\Omega}u_t(x,\tau)\mathrm{div}\hat{T}_E(x,\tau))dxd\tau\\[1ex]
&+&\int\nolimits_0^t\int\nolimits_{\Gamma_D}g_{D,t}(x,\tau)(T(x,\tau)-\hat{T}_E(x,\tau))\cdot n(x)dSd\tau+ \int\nolimits_0^t\int\nolimits_{\Omega}\D^{-1}T_{E,t}(x,\tau)\hat{T}_E(x,\tau)dxd\tau\nn\\[1ex]
&+& \frac{1}{c}\int\nolimits_0^t\int\nolimits_{\Omega}b_t(x,\tau)\hat{b}(x,\tau)dxd\tau 
+4l_c\int\nolimits_0^t\int\nolimits_{\partial\Omega}\nabla\mathrm{axl}(A(x,\tau))\cdot n \;\mathrm{axl}(A_{D,t}(x,\tau))dSd\tau\nn
\end{eqnarray}
is satisfied for all $t\in(0,T)$, where $T^0_E\in L^2(\Omega;\S)$ and $(u(0),A(0))\in H^{1}(\Omega;\R^3)\times H^2(\Omega;\mathfrak{so}(3))$ are unique solution of the problem (\ref{eq:3.8}).
\end{de}
\begin{tw}
$\mathrm{(Main\; existence\; result)}$\\
\label{tw:3.3}
Let us assume that the given data and initial data satisfy the properties, which are specified in (\ref{eq:3.5}) - (\ref{eq:3.8}). Then there exists a global in time solution (in the sense of Definition \ref{de:3.2}) of the system (\ref{eq:2.1})
with boundary conditions (\ref{eq:2.2}), (\ref{eq:2.3}) and initial condition (\ref{eq:2.4}).
\end{tw}
Notice that the solution defined above has a quite nice regularity. We even get that $\ve^p_t\in L^2(L^2)$, which yields $|\ve^p_t|b\in L^2(L^2)$. Unfortunately, this information is still not enough to obtain $L^2$-strong solutions. This paper presents the first existence result for the Armstrong-Frederick model with Cosserat effects and the new solution concept.\\
The proof of Theorem \ref{tw:3.3} is divided into two sections. First, we use the Yosida Approximation to the maximal monotone part of the inelastic constitutive equation. Next, we pass to the limit to obtain a solution in the sense of Definition \ref{de:3.2}.   
\section{Existence for the Yosida approximation}
We use the Yosida approximation for the monotone part of the flow rule from (\ref{eq:2.1}) in order to get a Lipschitz-nonlinearity only in equation $(\ref{eq:2.1})_4$. Hence, we obtain the following system of equations
\renewcommand{\theequation}{\thesection.\arabic{equation}}
\setcounter{equation}{0}%
\begin{eqnarray}
\label{eq:4.1}
\mathrm{div}_x T^{\nu}&=&-f\,,\nn\\
T^{\nu}&=&2\mu(\ve(u^{\nu})-\ve^{p,\nu})+2\mu_c(\mathrm{skew}(\nabla_x u^{\nu})-A^{\nu}) +\lambda\mathrm{tr}(\ve(u^{\nu})-\ve^{p,\nu})\id\,,\nn\\[1ex]
-l_c\,\Delta_x\mathrm{axl}\,(A^{\nu})&=&\mu_c\,\mathrm{axl}\,(\mathrm{skew}(\nabla_x u^{\nu})-A^{\nu})\,,\nn\\[1ex]
\ve^{p,\nu}_{t}&=&\frac{1}{\nu}\,\{|\dev\,(T_E^{\nu})-b^{\nu}|-\KK\}_{+}\frac{\dev\,(T_E^{\nu})-b^{\nu}}{|\dev\,(T_E^{\nu}) -b^{\nu}|}\,,\\[1ex]
T_E^{\nu}&=&2\mu(\ve(u^{\nu})-\ve^{p,\nu})+\lambda\mathrm{tr}(\ve(u^{\nu})-\ve^{p,\nu})\id,\nn\\[1ex]
b_t^{\nu}&=&c\,\ve^{p,\nu}_t-d\,|\ve_t^{p,\nu}|b^{\nu}.\nn
\end{eqnarray}
The above equations are studied for $x\in\Omega\subset \R^3$ and $t\in (0,T)$. $\nu>0$ and $\{\rho\}_{+}=\max\{0,\rho\}$, where $\rho$ is a 
scalar function.\\
The system (\ref{eq:4.1}) is considered with boundary conditions:
\begin{eqnarray}
\label{eq:4.2}
u^{\nu}(x,t)&=&g_D(x,t)\qquad \textrm{ for}\quad x\in\Gamma_D \quad\textrm{and}\quad t\geq 0\textrm{,}\nn\\
T^{\nu}(x,t)\cdot n(x)&=&g_N(x,t)\qquad\; \textrm{for}\quad x\in\Gamma_N \quad\textrm{and}\quad t\geq 0\textrm{,}\nn\\
A^{\nu}(x,t)&=&A_D(x,t)\quad\;\;\; \textrm{ for}\quad x\in \partial \Omega \quad\textrm{and}\quad  t\geq 0
\end{eqnarray}
and initial conditions
\begin{eqnarray}
\label{eq:4.3}
\ve^{p,\nu}(x,0)=\ve^{p,0}(x),\qquad b^{\nu}(x,0)=b^0(x).
\end{eqnarray}
\begin{tw}
\label{tw:4.1}
Fix $T>0$. Suppose that all hypotheses of Theorem \ref{tw:3.3} are satisfied. Then for all $\nu>0$ there exists a unique $L^2$- strong solution (in the sense of Definition \ref{de:3.1})  
$$(u^{\nu},T^{\nu},A^{\nu},\ve^{p,\nu},b^{\nu})\in W^{1,\infty}\Big(0,T;H^1(\Omega;\R^3)\times L^2(\Omega;\R^9)\times H^2(\Omega;\mathfrak{so}(3))\times (L^2(\Omega;\SS))^2\Big)$$
 satisfying the system (\ref{eq:4.1}) with boundary conditions (\ref{eq:4.2}) and initial conditions (\ref{eq:4.3}).   
\end{tw}
{\bf\large\em Proof:}\hspace{2ex} The proof is very similar to that of a related result of the Armstrong-Frederick model without Cosserat effects: see Section 4 of \cite{7}. It uses similar techniques. We provide a sketch below.\\[1ex]
{\bf Step 1:}\hspace{2ex} In the beginning we formulate two lemmas.
\begin{lem}
\label{lem:4.2}
Fix $T>0$. Assume that $(u^{\nu},T^{\nu},A^{\nu},\ve^{p,\nu},b^{\nu})$ is $L^2$- strong solution of the problem (\ref{eq:4.1}) and $|b^0(x)|\leq \frac{c}{d}$ for almost all $x\in\Omega$. Then for all $\nu>0$
$$|b^{\nu}(x,t)|\leq \frac{c}{d}\quad\mathrm{for\; a.\;e.\;}\quad (x,t)\in\Omega\times (0,T).$$ 
\end{lem}
The Lemma \ref{lem:4.2} implies that if $|b^0(x)|\leq \frac{c}{d}$ then we can modify the equation $(\ref{eq:4.1})_6$ in the following form 
\begin{equation}
\label{eq:4.4}
b_t^{\nu}=c\,\ve^{p,\nu}_t-d\,|\ve_t^{p,\nu}|\Pi(b^{\nu}),
\end{equation}
where
\begin{displaymath}
\Pi(b)= \left\{ \begin{array}{ll}
b & \mathrm{if}\quad |b|\leq \frac{c}{d}\,, \\[1ex]
\frac{c}{d}\,\frac{b}{|b|}& \mathrm{if}\quad |b|> \frac{c}{d}\,.
\end{array} \right.
\end{displaymath}
$\Pi$ is a projection on the convex set $\{b\in\SS:\; |b|\leq\frac{c}{d}\}$.
\begin{lem}
\label{lem:4.3}
Fix $T>0$. Assume that $(u^{\nu},T^{\nu},A^{\nu},\ve^{p,\nu},b^{\nu})$ satisfies the equation (\ref{eq:4.4}) in the $L^2$- strong sense on $\Omega\times (0,T)$ and $|b^0(x)|\leq \frac{c}{d}$ for almost all $x\in\Omega$. Then for all $\nu>0$
$$\Pi(b^{\nu}(x,t))=b^{\nu}(x,t).$$ 
\end{lem}
For the proofs of Lemma \ref{lem:4.2} and \ref{lem:4.3} we refer to Section 4, Lemma 1 and 2 in \cite{7}. From the above lemmas we conclude that if we find $L^2$-strong solution for the system (\ref{eq:4.1}) with the equation (\ref{eq:4.4}) instead of $(\ref{eq:4.1})_6$ then this solution will be also solution of the original system (\ref{eq:4.1}). Using this information we define the function $G_{\nu}:\S\times\SS\rightarrow \SS\times\SS$ by the formula  
\begin{eqnarray}
\label{eq:4.5}
G_{\nu}(T,b)= \left\{ \begin{array}{l}
\frac{1}{\nu}\,\{|\dev\,(T)-c\,b|-\KK\}_{+}\frac{\dev\,(T)-c\,b}{|\dev\,(T) -c\,b|}\,,\\[3ex]
\frac{c}{\nu}\,\{|\dev\,(T)-c\,b|-\KK\}_{+}\frac{\dev\,(T)-c\,b}{|\dev\,(T)-c\,b|}+\frac{cd}{\nu}\,\{|\dev\,(T)-c\,b|-\KK\}_{+}\Pi(b)\,.
\end{array} \right.
\end{eqnarray}
Let us consider the system (\ref{eq:4.1}) with the equation (\ref{eq:4.4}) instead of $(\ref{eq:4.1})_6$. The flow rule for the modified system can by written in the form
\begin{eqnarray}
\label{eq:4.6}
(\ve^{p,\nu}_t,b^{\nu}_t)=G_{\nu}(-\rho\nabla_{(\ve^{p},b)}\psi^{\nu}(\ve^{\nu},\ve^{p,\nu},A^{\nu},b^{\nu})),
\end{eqnarray}
where the free energy function $\psi^{\nu}$ is in the form
\begin{eqnarray*}
\rho\psi^{\nu}(\ve^{\nu},\ve^{p,\nu},A^{\nu},b^{\nu})&=&\mu\|\ve(u^{\nu})-\ve^{p,\nu}\|^2+
\mu_c\|\mathrm{skew}(\nabla_x u^{\nu})-A^{\nu}\|^2\nn\\
&+&\frac{\lambda}{2}\Big(\mathrm{tr}(\ve(u^{\nu})-\ve^{p,\nu})\Big)^2+2l_c\|\nabla_x\mathrm{axl}(A^{\nu})\|^2+\frac{1}{2c}\|b^{\nu}\|^2\,.
\end{eqnarray*}
Denote by $\E^{\nu}(t)$ the energy associated with the system (\ref{eq:4.1})
\begin{eqnarray}
\label{eq:4.20}
\E^{\nu}(u^{\nu}\ve^{\nu},\ve^{p,\nu},A^{\nu},b^{\nu})(t)=\int\nolimits_{\Omega}
\rho\psi^{\nu}\Big(u^{\nu}(x,t),\ve^{\nu}(x,t),\ve^{p,\nu}(x,t),A^{\nu}(x,t),b^{\nu}(x,t)\Big)dx.
\end{eqnarray}
A fundamental tool in our proof is the following property of the energy function which results from our Cosserat modification:
\begin{tw}
\label{tw:4.7}
$\mathrm{(coerciveness\; of\; the\; energy)}$\\
(a) $\mathit{(the\; case\; with\; zero\; boundary\; data)}$\\
For all $\nu>0$ the energy function (\ref{eq:4.20}) is elastically coercive with respect to $\nabla u$. 
This means that $\exists\; C_E>0$, $\forall\; u\in H^1_0(\Omega)$, $\forall\; A\in H^1_0(\Omega)$, $\forall\; \ve^p\in L^2(\Omega)$, $\forall\; b\in L^2(\Omega)$
$$\E^{\nu}(u,\ve,\ve^{p},A,b)\geq C_E\Big(\|u\|^2_{H^1(\Omega)}+\|A\|^2_{H^1(\Omega)}+\|b\|^2_{L^2(\Omega)}\Big).$$ 
(b) $\mathit{(the\; case\; with\; non-zero\; boundary\; data)}$\\
Moreover, $\exists\; C_E>0$, $\forall\; g_D,\;A_D\in H^{\frac{1}{2}}(\partial\Omega)$, 
$\exists\; C_D>0$, $\forall\; \ve^p\in L^2(\Omega)$, $\forall\; b\in L^2(\Omega)$, $\forall\; u\in H^1(\Omega)$, $\forall\; A\in H^1(\Omega)$ 
with boundary conditions $u_{|_{\partial\Omega}}=g_D$ and $A_{|_{\partial\Omega}}=A_D$ it holds that
$$\E^{\nu}(u,\ve,\ve^{p},A,b)+C_D\geq C_E\Big(\|u\|^2_{H^1(\Omega)}+\|A\|^2_{H^1(\Omega)}+\|b\|^2_{L^2(\Omega)}\Big).$$
\end{tw}
For the proof of Theorem \ref{tw:4.7} we refer to the Theorem $3.2$ of the article \cite{8}. Moreover,
Lemma 3 of \cite{7} provides the following properties of $G_{\nu}$:\\[1ex]
\textbullet\hspace{1ex} {\em The function $G_{\nu}$ generates a bounded nonlinear operator from $L^2(\Omega;\S\times\SS)$ into $L^2(\Omega;\S\times\SS)$.}\\[1ex]
\textbullet\hspace{1ex} {\em For all $(T^{1},b^1)$, $(T^{2},b^2)\in L^2(\Omega;\S\times\SS)$ the inequality
\begin{eqnarray*}
\Big(G_{\nu}(T^{1},b^1)-G_{\nu}(T^{2},b^2),(T^{1},b^1)-(T^{2},b^2)\Big)_{L^2} \geq -\frac{c^2\,(c+\frac{1}{2})}{\nu}\Big(\|T^1-T^2\|^2_{L^2}+\|b^1-b^2\|^2_{L^2}\Big)
\end{eqnarray*}
holds.}\\[1ex]
By the first statement above we can conclude that the system (\ref{eq:4.1}) with the flow rule (\ref{eq:4.5}) possesses the linear self-controlling property (for the definition we refer to \cite{6}) and the second one means that this flow rule belongs to the class $\cal L\cal M$ - the class of Lipschitz perturbations of monotone vector fields.\\[1ex]
{\bf Step 2:}\hspace{2ex} The reasoning from above gives us that the system (\ref{eq:4.1}) with the flow rule (\ref{eq:4.5}) is coercive in the sense of Theorem \ref{tw:4.7} and it belongs to the class $\cal L\cal M$ with self-controlling property. The statement of the Theorem \ref{tw:4.1} may be proved using the general theory developed for the class $\cal L\cal M$ with self-controlling property (see for example \cite{5,6}). This theory was also used to study the existence of solutions to a model of poroplasticity with Cosserat effect, which has a similar structure as the system (\ref{eq:2.1}) (details can be found in Section 4 of \cite{11}). However, we would like very briefly to present the theory, which was used for the Armstrong-Frederick model without Cosserat effects (we refer to \cite{7} for details). Let $k>0$ be a positive real number and $C^k$ be the following cut function 
\begin{displaymath}
C^k:\R_+\rightarrow\R_+\;,\quad C^k(s)= \left\{ \begin{array}{ll}
s & \mathrm{if}\quad s\leq k\,, \\
k & \mathrm{if}\quad s> k\,.
\end{array} \right.
\end{displaymath}
Consider the following sequence of problems
\begin{eqnarray}
\label{eq:4.7}
\mathrm{div}_x T^{k}&=&-f\,,\nn\\
T^{k}&=&2\mu(\ve(u^{k})-\ve^{p,k})+2\mu_c(\mathrm{skew}(\nabla_x u^{k})-A^{k}) +\lambda\mathrm{tr}(\ve(u^{k})-\ve^{p,k})\,\id\,,\nn\\[1ex]
-l_c\,\Delta_x\mathrm{axl}\,(A^{k})&=&\mu_c\,\mathrm{axl}\,(\mathrm{skew}(\nabla_x u^{k})-A^{k})\,,\nn\\[1ex]
\ve^{p,k}_{t}&=&\frac{1}{\nu}\,\{|\dev\,(T_E^{k})-b^{k}|-\KK\}_{+}\frac{\dev\,(T_E^{k})-b^{k}}{|\dev\,(T_E^{k}) -b^{k}|}\,,\\[1ex]
T_E^{k}&=&2\mu(\ve(u^{k})-\ve^{p,k})+\lambda\mathrm{tr}(\ve(u^{k})-\ve^{p,k})\,\id\,,\nn\\[1ex]
b_t^{k}&=&c\,C^k(|\ve^{p,k}_{t}|)\frac{\ve^{p,k}_t}{|\ve^{p,k}_{t}|}-d\,C^k(|\ve_t^{p,k}|)\Pi(b^{k})\,.\nn
\end{eqnarray}
Here we drop the superscript $\nu>0$ and write $(u^{k},A^{k},T^{k},\ve^{p,k},b^{k})$  instead of\\ $(u^{\nu,k},A^{\nu,k},T^{\nu,k},\ve^{p,\nu,k},b^{\nu,k})$. The system (\ref{eq:4.7}) is considered with boundary conditions (\ref{eq:4.2}) and initial condition (\ref{eq:4.3}).
\begin{tw}
\label{tw:4.4}
$\mathrm{(global\; existence\; for\; Lipschitz\; nonlinearities)}$\\
Assume that for all $T>0$ the given data has the following regularity
\begin{eqnarray*}
f\in C^{1}((0,T];L^2(\Omega;\R^3)), && g_D\in C^{1}([0,T];H^{\frac{1}{2}}(\Gamma_D;\R^3))\,,\\[1ex]
g_N\in C^{1}([0,T];H^{-\frac{1}{2}}(\Gamma_N;\R^3)),&&
A_D\in C^{1}([0,T];H^{\frac{3}{2}}(\partial\Omega;\so(3)))\,,\\[1ex]
(\ve^{p,0},b^0)&\in& L^2(\Omega;\SS)\times L^2(\Omega;\SS)\,. 
\end{eqnarray*}
Then, for all $k>0$ the approximate problem (\ref{eq:4.7}) with boundary conditions (\ref{eq:4.2}) and initial condition (\ref{eq:4.3}) has a global in time, unique $L^2$ - strong solution $(u^{k},A^{k},T^{k},\ve^{p,k},b^{k})$ with regularity
\begin{eqnarray*}
(u^{k},A^{k},T^k)\in C^{1}([0,T];H^1(\Omega;\R^3)\times H^2(\Omega;\so(3))\times L^2(\Omega;\R^9)),\nn
\end{eqnarray*}
\begin{eqnarray*}
\ve^{p,k}\in C^1([0,T]; L^2(\Omega;\SS)),\quad b^{k}\in C^1([0,T]; L^2(\Omega;\SS).
\end{eqnarray*}
\end{tw}
Notice that the function $C^k(|\ve^{p,k}_{t}|)$ is an $L^{\infty}$-function therefore, the system (\ref{eq:4.7}) includes global Lipschitz nonlinearities only and for a proof of this statement we refer to Theorem 3.1 of the article \cite{8}.\\[1ex]
{\bf Step 3:}\hspace{2ex} Our goal is to pass to the limit with $k\rightarrow\infty$. To do this we use the energy method and the properties of the function $G_{\nu}$. 
\begin{tw}
\label{tw:4.5}
Let us suppose that all hypotheses of Theorem \ref{tw:4.1} hold. Let us denote by $\E^{k}(t)$ the total energy associated with the system (\ref{eq:4.7})
\begin{eqnarray*}
\E^{k}(\ve^{k},\ve^{p,k},A^{k},b^{k})(t)&=&\int\nolimits_{\Omega}\Big( \mu|\ve(u^k)(x,t)-\ve^{p,k}(x,t)|^2+
\mu_c|\mathrm{skew}(\nabla_x u^k(x,t))-A^k(x,t)|^2\nn\\
&+&\frac{\lambda}{2}\Big(\mathrm{tr}(\ve(u^k(x,t))-\ve^{p,k}(x,t))\Big)^2+2l_c|\nabla_x\mathrm{axl}(A^k(x,t))|^2\\
&+&\frac{1}{2c}|b^k(x,t)|^2\Big)dx.
\end{eqnarray*}
Then for all $t\in(0,T)$ the following estimate
\begin{eqnarray*}
\E^{k}(\ve^{k},\ve^{p,k},A^{k},b^{k})(t)\leq C(T)
\end{eqnarray*}
holds and $C(T)$ does not depend on $k>0$ (it depends only on the given data and the domain).
\end{tw}
{\bf\em Proof of Theorem \ref{tw:4.5}:}\hspace{2ex}Calculate the time derivative of the energy and obtain
\begin{eqnarray}
\label{eq:4.8}
&&\frac{d}{dt}\Big(\E^{k}(\ve^{k},\ve^{p,k},A^{k},b^k)(t)\Big)=
\int\nolimits_{\Omega}\Big(2\mu(\ve^{k}-\ve^{p,k})(\ve^{k}_t-\ve_t^{p,k})
+\lambda\mathrm{tr}(\ve^{k}-\ve^{p,k})\mathrm{tr}(\ve_t^{k}-\ve_t^{p,k})\nn\\[1ex]
&+&2\mu_c(\mathrm{skew}(\nabla u^{k})-A^{k})(\mathrm{skew}(\nabla u_t^{k})-A_t^{k})\Big)\,dx
+4l_c\int\nolimits_{\Omega}\nabla\mathrm{axl}(A^{k})\nabla\mathrm{axl}(A_t^{k})\,dx\nn\\[1ex]
&+&\frac{1}{c}\int\nolimits_{\Omega}b^kb^k_tdx =\int\nolimits_{\Omega}T^{k}\nabla u^{k}_t\,dx
-\int\nolimits_{\Omega}\dev\,(T_E^{k})\ve_t^{p,k}\,dx-2\mu_c\int\nolimits_{\Omega}(\mathrm{skew}(\nabla u^{k})-A^{k})A_t^{k}\,dx\nn\\[1ex]
&+&4l_c\int\nolimits_{\Omega}\nabla\mathrm{axl}(A^{k})\nabla\mathrm{axl}(A_t^{k})\,dx -
\int\nolimits_{\Omega}b^k\Big(\ve_t^{p,k}-C^k(|\ve_t^{p,k}|)\frac{\ve_t^{p,k}}{|\ve_t^{p,k}|}\Big)dx\nn\\[1ex]
&&\underbrace{-\int\nolimits_{\Omega}|\dev\,(T_E^{k})-b^k|\,|\ve_t^{p,k}|\,dx - 
\frac{c}{d}\int\nolimits_{\Omega}C^k(|\ve_t^{p,k}|)|\Pi(b^k)|\,|b^k|dx}_{\leq 0}.
\end{eqnarray} 
Integrating by parts in the first and fourth term on the right hand side of (\ref{eq:4.8}), using the equations $(\ref{eq:4.7})_1$ and $(\ref{eq:4.7})_3$ and boundary data we have (notice that $\|A\|^2=2\|\mathrm{axl}A\|^2$)
\begin{eqnarray}
\label{eq:4.9}
&&\frac{d}{dt}\Big(\E^{k}(\ve^{k},\ve^{p,k},A^{k},b^k)(t)\Big)\leq
\int\nolimits_{\Omega}fu^{k}_t\,dx +
\int\nolimits_{\Gamma_D}g_{D,t}T^{k}\cdot ndS+\int\nolimits_{\Gamma_N}u^{k}_tg_{N}dS\nn\\[1ex]
&+&4l_c\int\nolimits_{\partial\Omega}\nabla\mathrm{axl}(A^{k})\cdot n \;\mathrm{axl}(A_{D,t})dS
-\int\nolimits_{\Omega}b^k\Big(\ve_t^{p,k}-C^k(|\ve_t^{p,k}|)\frac{\ve_t^{p,k}}{|\ve_t^{p,k}|}\Big)dx.
\end{eqnarray}
Integrating (\ref{eq:4.9}) with respect to time we have
\begin{eqnarray}
\label{eq:4.10}
&& \E^{k}(\ve^{k},\ve^{p,k},A^{k},b^k)(t)\leq\qquad
\E^{k}(\ve^{k},\ve^{p,k},A^{k},b^k)(0)\nn\\[1ex]
&+&\int\nolimits_0^t\int\nolimits_{\Omega}f(\tau)u^{k}_t(\tau)\,dxd\tau +
\int\nolimits_0^t\int\nolimits_{\Gamma_D}g_{D,t}(\tau)T^{k}(\tau)\cdot ndSd\tau\nn\\[1ex]
&+&\int\nolimits_0^t\int\nolimits_{\Gamma_N}u^{k}_t(\tau)g_{N}(\tau)dS
+4l_c\int\nolimits_0^t\int\nolimits_{\partial\Omega}\nabla\mathrm{axl}(A^{k}(\tau))\cdot n \;\mathrm{axl}(A_{D,t}(\tau))dS\nn\\[1ex]
&-&\int\nolimits_0^t\int\nolimits_{\Omega}b^k(\tau)\Big(\ve_t^{p,k}(\tau)-C^k(|\ve_t^{p,k}(\tau)|)\frac{\ve_t^{p,k}(\tau)}{|\ve_t^{p,k}(\tau)|}\Big)dxd\tau.
\end{eqnarray}
The continuity with respect to time yields that the initial values $u^k(0)$, $A^k(0)$ are solutions of the following linear elliptic boundary-value problem
\begin{eqnarray*}
\mathrm{div}_x T^k(0)&=&-f(0)\,,\nn\\[1ex]
-l_c\,\Delta_x\mathrm{axl}\,(A^k(0))&=&\mu_c\,\mathrm{axl}\,(\mathrm{skew}(\nabla_x u^k(0))-A^k(0))\,,\nn\\[1ex]
u^k(0)_{|_{\Gamma_D}}=g_D(0) && T^k(x)\cdot n_{|_{\Gamma_N}}=g_N(0)\quad A^k(0)_{|_{\partial\Omega}}=A_D(0)\,,
\end{eqnarray*}
where 
$$T^k(0)=2\mu(\ve(u^k(0))-\ve^{p,0})+2\mu_c(\mathrm{skew}(\nabla_x u^k(0))-A^k(0))+\lambda\mathrm{tr}(\ve(u^k(0))-\ve^{p,0})\id.\nn\\$$
This solution has the following regularity
$$u^k(0)\in H^1(\Omega;\R^3),\qquad A^k(0)\in H^2(\Omega;\so(3))$$
and is independent on $k$. Hence, the initial energy $\E^{k}(\ve^{k},\ve^{p,k},A^{k},b^k)(0)$ is a constant. Integrating partially in time, we get
\begin{eqnarray}
\label{eq:4.11}
&&\int\nolimits_0^t\int\nolimits_{\Omega}f(\tau)u^{k}_t(\tau)\,dxd\tau= -\int\nolimits_0^t\int\nolimits_{\Omega}f_t(\tau)u^{k}(\tau)\,dxd\tau+ \int\nolimits_{\Omega}f(t)u^{k}(t)\,dx+\int\nolimits_{\Omega}f(0)u^{k}(0)\,dx\nn\\[1ex]
&\leq&\frac{1}{2}\int\nolimits_0^t\|f_t(\tau)\|^2_{L^2(\Omega;\R^3)}d\tau+\frac{1}{2}\int\nolimits_0^t\|u^{k}(\tau)\|^2_{L^2(\Omega;\R^3)}d\tau+ \|f(t)\|_{L^2(\Omega;\R^3)}\|u^{k}(t)\|_{L^2(\Omega;\R^3)}\nn\\[1ex]
&+&\|f(0)\|_{L^2(\Omega;\R^3)}\|u^{k}(0)\|_{L^2(\Omega;\R^3)}\,.
\end{eqnarray}
Applying Poincar\'e's inequality to $\|u^k\|_{L^2}$ in (\ref{eq:4.11}) and the coerciveness of the energy function with respect to the gradient of the displacement vector we obtain
\begin{eqnarray}
\label{eq:4.12}
&&\Big|\int\nolimits_0^t\int\nolimits_{\Omega}f(\tau)u^{k}_t(\tau)\,dxd\tau\Big|\leq\qquad \hat{C}\int\nolimits_0^t\E^{k}(\ve^{k},\ve^{p,k},A^{k},b^k)(\tau)d\tau\nn\\[1ex]
&&\qquad\qquad+\hat{C}\|f(t)\|_{L^2(\Omega;\R^3)}\Big(\E^{k}(\ve^{k},\ve^{p,k},A^{k},b^k)(t)\Big)^{\frac{1}{2}}+\hat{C}(t)\,,
\end{eqnarray}
where the constants $\hat{C}$, $\hat{C}(t)$ do not depend on $k$. In a similar way we show that
\begin{eqnarray}
\label{eq:4.13}
&&\Big|\int\nolimits_0^t\int\nolimits_{\Omega}g_{N}(\tau)u^{k}_t(\tau)\,dSd\tau\Big|\leq\qquad \hat{C}\int\nolimits_0^t\E^{k}(\ve^{k},\ve^{p,k},A^{k},b^k)(\tau)d\tau\nn\\[1ex]
&&\qquad\qquad+\hat{C}\|g_N(t)\|_{H^{-\frac{1}{2}}(\Gamma_N;\R^3)}\Big(\E^{k}(\ve^{k},\ve^{p,k},A^{k},b^k)(t)\Big)^{\frac{1}{2}}+\hat{C}(t)\,.
\end{eqnarray}
The others boundary integrals in (\ref{eq:4.10}) are estimated using the trace theorems (the details will be shown in the proof of Theorem \ref{tw:5.1} of the next Section). It is not difficult to see that Lemma \ref{lem:4.2} holds for the sequence $b^k$ and 
\begin{eqnarray}
\label{eq:4.14}
-\int\nolimits_{\Omega}b^k\Big(\ve_t^{p,k}-C^k(|\ve_t^{p,k}|)\frac{\ve_t^{p,k}}{|\ve_t^{p,k}|}\Big)dx&\leq& \int\nolimits_{\Omega}|b^k||\ve_t^{p,k}|dx\leq \|b^k\|_{L^{\infty}(\Omega;\SS)} \|\ve^{p,k}_t\|_{L^{1}(\Omega;\SS)}\nn\\[1ex]
&\leq&C(1+\alpha\|\ve^{p,k}_t\|^2_{L^{2}(\Omega;\SS)})\nn\\[1ex]
&\leq& \quad(\quad\mathrm{equation\;\;\;}(4.1)_4\quad)\quad\leq\\[1ex]
&\leq&C\Big(1+\frac{\alpha}{\nu}(\|T^k\|^2_{L^{2}(\Omega;\S)}+\|b^{k}\|^2_{L^{2}(\Omega;\SS)})\Big)\nn\,,
\end{eqnarray}
where $\alpha$ is any positive number and $C$ does not depend on $k$. By (\ref{eq:4.11})-(\ref{eq:4.14}) we are able to obtain the following inequality 
\begin{eqnarray}
\label{eq:4.15}
\E^{k}(\ve^{k},\ve^{p,k},A^{k},b^k)(t)&\leq&\qquad \alpha\,\E^{k}(\ve^{k},\ve^{p,k},A^{k},b^k)(t)\nn\\[1ex]
&+&C\int\nolimits_0^t \E^{k}(\ve^{k},\ve^{p,k},A^{k},b^k)(\tau)d\tau+C(T)\,,
\end{eqnarray}
where $\alpha$ is any positive number and the constant $C$, $C(T)>0$ does not depend on $k>0$. Choosing $\alpha>0$ sufficiently small and using Gronwall's lemma we complete the proof of Theorem \ref{tw:4.5}.$\mbox{}$ \hfill $\Box$\\[1ex]
From the proof of the Theorem \ref{tw:4.5} we conclude that the sequence $(u^{k},A^{k},T^{k},\ve^{p,k},b^{k})$ is bounded in $W^{1,\infty}(H^1\times H^2\times L^2\times L^2 \times L^2)$.\\[1ex]
{\bf Step 4:}\hspace{2ex} To finish the proof of Theorem \ref{tw:4.1} we need $L^{\infty}(L^2)$-strong convergence of the sequence $(T^{k},b^{k})$.
\begin{tw}
\label{tw:4.6}
Let us assume that the given data satisfy all requirements of Theorem \ref{tw:4.5}. Then,
\begin{eqnarray*}
\E^{k}(\ve^{m}-\ve^{n},\ve^{p,m}-\ve^{p,n},A^{m}-A^{n},b^{m}-b^{n})(t) \longrightarrow  0
\end{eqnarray*}
for $m$,$n\rightarrow\infty$ uniformly on bounded time intervals. 
\end{tw}
For the proof of the above Theorem we refer to Section $4$ of \cite{7}. Theorem \ref{tw:4.6} implies the $L^{\infty}(L^2)$-strong convergence of the sequence $(T^{k},b^{k})$ hence we can pass to the limit in the system (\ref{eq:4.7}) with $k\rightarrow\infty$ and obtain $L^2$-strong solution of the system (\ref{eq:4.1}). The uniqueness follows immediately from coerciveness of the energy function evaluated on the difference of two solutions of the system (\ref{eq:4.1}). The last statement finishes the proof of Theorem \ref{tw:4.1}.$\mbox{}$ \hfill $\Box$\
\section{Proof of the Main Theorem \ref{tw:3.3}}
In this section we are going to prove the $L^2(L^2)$- boundedness for the time derivatives of the sequence $(u^{\nu},T^{\nu},A^{\nu},\ve^{p,\nu},b^{\nu})$. It is the main part to prove
 the Theorem \ref{tw:3.3}.
\begin{tw}$\mathrm{(Energy\; estimate)}$\\
\label{tw:5.1}
Assume that the given data and initial data satisfies (\ref{eq:3.5}) - (\ref{eq:3.8}). Then for all $t\in(0,T)$ the following estimate
\begin{eqnarray*}
&&\int\nolimits_{\Omega}\frac{1}{2\nu}\{|\dev\,(T_E^{\nu})(t)-b^{\nu}(t)|-\KK\}_{+}^2dx+\int\nolimits^t_0\int\nolimits_{\Omega}\D^{-1}T_{E,t}^{\nu}(\tau)T_{E,t}^{\nu}(\tau)dxd\tau\\[1ex]
&+& 2\mu_c\int\nolimits^t_0\int\nolimits_{\Omega}|\mathrm{skew}(\nabla_x u_t^{\nu}(\tau))-A_t^{\nu}(\tau)|^2dxd\tau
+4l_c\int\nolimits^t_0\int\nolimits_{\Omega}|\nabla\mathrm{axl}(A_t^{\nu}(\tau))|^2dxd\tau
\leq C(T)
\end{eqnarray*}
holds and $C(T)$ does not depend on $\nu>0$ (it depends only on the given data and the domain).
\end{tw}
{\bf\em Proof:}\hspace{2ex} Compute the time derivative 
\renewcommand{\theequation}{\thesection.\arabic{equation}}
\setcounter{equation}{0}%
\begin{eqnarray} 
\label{eq:5.1} 
&&\frac{d}{dt}\Big (\int\nolimits_{\Omega}\frac{1}{2\nu}\{|\dev\,(T_E^{\nu})(t)-b^{\nu}(t)|-\KK\}_{+}^2dx\Big )= 
\int\nolimits_{\Omega} \ve^{p,\nu}_t(t)\Big(\dev\,(T_{E,t}^{\nu})(t)-b^{\nu}_t(t)\Big)dx\nn\\[1ex]
&=&(\textrm{by}\;\textrm{the}\;\textrm{elastic}\;\textrm{constitutive}\;\textrm{relation})= \int\nolimits_{\Omega} \ve^{\nu}_t(t)T_{E,t}^{\nu}(t)dx-
\int\nolimits_{\Omega} \D^{-1}T_{E,t}^{\nu}T_{E,t}^{\nu}(t)dx\nn\\[1ex]
&-&\int\nolimits_{\Omega} \ve^{p,\nu}_t(t)b^{\nu}_t(t)dx.
\end{eqnarray}
From Lemma \ref{lem:4.2} we have the following inequality
\begin{eqnarray}
\label{eq:5.2}  
&&\int\nolimits_{\Omega} \ve^{p,\nu}_t(t)b^{\nu}_t(t)dx= 
(\textrm{by}\;\textrm{the}\;\textrm{equation}\;\textrm{for}\;\textrm{the}\;\textrm{backstress})=\\[1ex]
&&c\int\nolimits_{\Omega} |\ve^{p,\nu}_t(t)|^2dx
-d\int\nolimits_{\Omega} |\ve^{p,\nu}_t(t)|\ve^{p,\nu}_t(t)b^{\nu}(t)dx\geq 
\int\nolimits_{\Omega} |\ve^{p,\nu}_t(t)|^2(c-d\,|b^{\nu}(t)|)dx\geq 0\,.\nn
\end{eqnarray}
Notice that
\begin{eqnarray} 
\label{eq:5.3} 
&&\int\nolimits_{\Omega} \ve^{\nu}_t(t)T_{E,t}^{\nu}(t)dx=\int\nolimits_{\Omega} \nabla u^{\nu}_t(t)T_{t}^{\nu}(t)dx-
\int\nolimits_{\Omega} 2\mu_c\Big(\mathrm{skew}(\nabla_x u_t^{\nu}(t))-A_t^{\nu}(t)\Big)\mathrm{skew}(\nabla_x u^{\nu}_t(t))dx\nn\\[1ex]
&=&\int\nolimits_{\Omega} \nabla u^{\nu}_t(t)T_{t}^{\nu}(t)dx-
2\mu_c\int\nolimits_{\Omega}|\mathrm{skew}(\nabla_x u_t^{\nu}(t))-A_t^{\nu}(t)|^2dx\nn\\[1ex]
&-&2\mu_c\int\nolimits_{\Omega}\Big(\mathrm{skew}(\nabla_x u_t^{\nu}(t))-A_t^{\nu}(t)\Big)A_t^{\nu}(t)dx.
\end{eqnarray} 
Using $(\ref{eq:4.1})_3$, integrating by parts in (\ref{eq:5.3}) we have
\begin{eqnarray} 
\label{eq:5.4} 
&&\int\nolimits_{\Omega} \ve^{\nu}_t(t)T_{E,t}^{\nu}(t)dx=\int\nolimits_{\Omega} u^{\nu}_t(t)f_t(t)dx +
\int\nolimits_{\Gamma_D}g_{D,t}(t)(T_t^{\nu}(t))\cdot n\;dS+\int\nolimits_{\Gamma_N}u^{\nu}_t(t)g_{N,t}(t)dS\nn\\[1ex]
&-&2\mu_c\int\nolimits_{\Omega}|\mathrm{skew}(\nabla_x u_t^{\nu}(t))-A_t^{\nu}(t)|^2dx- 
4l_c\int\nolimits_{\Omega}|\nabla\mathrm{axl}(A_t^{\nu}(t))|^2dx\nn\\[1ex]
&+&4l_c\int\nolimits_{\partial\Omega}\nabla\mathrm{axl}(A^{\nu}_t(t))\cdot n \;\mathrm{axl}(A_{D,t}(t))dS\,.
\end{eqnarray}
The first term on the right hand side of (\ref{eq:5.4}) is estimated as follows
\begin{eqnarray} 
\label{eq:5.19} 
\int\nolimits_{\Omega} u^{\nu}_t(t)f_t(t)dx \leq \|u_t^{\nu}(t)\|_{H^1(\Omega;\R^3)}\|f_t(t)\|_{L^2(\Omega;\R^3)}\,.
\end{eqnarray}
We estimate the appearing boundary integrals
\begin{eqnarray}
\label{eq:5.5}
&&\int\nolimits_{\Gamma_D}g_{D,t}(t)(T_t^{\nu}(t))\cdot ndS\leq \|g_{D,t}(t)\|_{H^{\frac{1}{2}}(\Gamma_D;\R^3)}
\|T_t^{\nu}(t)\cdot n\|_{H^{-\frac{1}{2}}(\Gamma_D;\R^3)}\nn\\[2ex]
&\leq&\quad(\textrm{from the trace theorem in the space }
L^2_{\mathrm{div}}(\Omega) \textrm{ see for example in \cite{14}})\quad\leq\nn\\[2ex]
&\leq&C\Big(\|T_t^{\nu}(t)\|_{L^2(\Omega;\S)}+\|\textrm{div} T_t^{\nu}(t)\|_{L^2(\Omega;\R^3)}\Big)
\|g_{D,t}(t)\|_{H^{\frac{1}{2}}(\partial\Omega;\R^3)}\,\nn\\
&\leq&C\|T_t^{\nu}(t)\|_{L^2(\Omega;\S)}\|g_{D,t}(t)\|_{H^{\frac{1}{2}}(\partial\Omega;\R^3)}+
C\|f_t(t)\|_{L^2(\Omega;\R^3)}
\|g_{D,t}(t)\|_{H^{\frac{1}{2}}(\partial\Omega;\R^3)}\,,
\end{eqnarray}
where $C>0$ does not depend on $\nu$. Moreover
\begin{eqnarray}
\label{eq:5.6}
\int\nolimits_{\Gamma_N}u^{\nu}_t(t)g_{N,t}(t)dS&\leq& \|g_{N,t}(t)\|_{H^{-\frac{1}{2}}(\Gamma_N;\R^3)}
\|u_t^{\nu}(t)\|_{H^{\frac{1}{2}}(\partial\Omega;\R^3)}\nn\\
&\leq&\|g_{N,t}(t)\|_{H^{-\frac{1}{2}}(\Gamma_N;\R^3)}
\|u_t^{\nu}(t)\|_{H^{1}(\Omega;\R^3)}\,.
\end{eqnarray}
We use the $H^2$-regularity of the microrotations to estimate the following integral
\begin{eqnarray}
\label{eq:5.7}
&&\int\nolimits_{\partial\Omega}(\nabla\mathrm{axl}(A^{\nu}_t(t)))\cdot n\,\mathrm{axl}(A_{D,t}(t))\,dS\nn\\[1ex]&\leq&\|\nabla\mathrm{axl}(A^{\nu}_t(t)))\cdot n\|_{H^{-\frac{1}{2}}(\partial\Omega;\R^3)}
\|\mathrm{axl}(A_{D,t}(t))\|_{H^{\frac{1}{2}}(\partial\Omega;\R^3)}\,\nn\\[2ex]
&\leq&C\Big(\|\nabla\mathrm{axl}(A^{\nu}_t(t))\|_{L^2(\Omega;\R^9)}
+ \|\Delta\mathrm{axl}(A^{\nu}_t(t))\|_{L^2(\Omega;\R^3)}\Big)\|\mathrm{axl}(A_{D,t}(t))\|_{H^{\frac{1}{2}}(\partial\Omega;\R^3)}\,\nn\\
&=&C\|\nabla\mathrm{axl}(A^{\nu}_t(t))\|_{L^2(\Omega;\R^9)}\|\mathrm{axl}(A_{D,t}(t))\|_{H^{\frac{1}{2}}(\partial\Omega;\R^3)}\nn\\
&+& C\frac{\mu_c}{l_c}\|\mathrm{skew}(\nabla u^{\nu}_t(t))-A^{\nu}_t(t)\|_{L^2(\Omega;\R^9)}
\|\mathrm{axl}(A_{D,t}(t))\|_{H^{\frac{1}{2}}(\partial\Omega;\R^3)}\,.
\end{eqnarray}
Inserting (\ref{eq:5.2})-(\ref{eq:5.7}) into (\ref{eq:5.1}), using Cauchy's inequality with a small weight at the approximate sequence and integrating with respect to time we obtain the following inequality
\begin{eqnarray}
\label{eq:5.8}
&&\int\nolimits_{\Omega}\frac{1}{2\nu}\{|\dev\,(T_E^{\nu})(t)-b^{\nu}(t)|-\KK\}_{+}^2dx+\int\nolimits^t_0\int\nolimits_{\Omega}\D^{-1}T_{E,t}^{\nu}(\tau)T_{E,t}^{\nu}(\tau)dxd\tau\nn\\[1ex]
&+& 2\mu_c\int\nolimits^t_0\int\nolimits_{\Omega}|\mathrm{skew}(\nabla_x u_t^{\nu}(\tau))-A_t^{\nu}(\tau)|^2dxd\tau
+4l_c\int\nolimits^t_0\int\nolimits_{\Omega}|\nabla\mathrm{axl}(A_t^{\nu}(\tau))|^2dxd\tau\nn\\[1ex]
&\leq & \alpha\int\nolimits^t_0\Big(\|T_t^{\nu}(\tau)\|^2_{L^2(\Omega;\S)} + \|\mathrm{skew}(\nabla u^{\nu}_t(\tau))-A^{\nu}_t(\tau)\|^2_{L^2(\Omega;\R^9)} + \|\nabla\mathrm{axl}(A^{\nu}_t(\tau))\|^2_{L^2(\Omega;\R^9)}\Big)d\tau\nn\\[1ex]
&+&\alpha\int\nolimits^t_0\|u_t^{\nu}(\tau)\|^2_{H^{1}(\Omega;\R^3)}d\tau+\tilde{C}(T,\alpha) + \int\nolimits_{\Omega}\frac{1}{2\nu}\{|\dev\,(T_E^{\nu})(0)-b^{\nu}(0)|-\KK\}_{+}^2dx\,,
\end{eqnarray}
where $\alpha>0$ is any positive number and $\tilde{C}(T)>0$ does not depend on $\nu$ (it depends on the given data and the time interval only). We observe that $T^{\nu}(0)\in L^2(\Omega;\S)$ is the unique solution of the problem
\begin{eqnarray*}
\mathrm{div}_x T^{\nu}(x,0)&=&-f(x,0)\,,\nn\\[1ex]
-l_c\,\Delta_x\mathrm{axl}\,(A^{\nu}(x,0))&=&\mu_c\,\mathrm{axl}\,(\mathrm{skew}(\nabla_x u^{\nu}(x,0))-A^{\nu}(x,0))\,,\nn\\[1ex]
u^{\nu}(x,0)_{|_{\Gamma_D}}=g_D(x,0)\,, && T^{\nu}(x,0)\cdot n_{|_{\Gamma_N}}=g_N(x,0),\quad A^{x,0}(x,0)_{|_{\partial\Omega}}=A_D(x,0)\,,
\end{eqnarray*}
where 
\begin{eqnarray*}
T^{\nu}(x,0)&=&2\mu(\ve(u^{\nu}(x,0))-\ve^{p,0}(x))+2\mu_c(\mathrm{skew}(\nabla_x u^{\nu}(x,0))-A^{\nu}(x,0))\\[1ex]
&+&\lambda\mathrm{tr}(\ve(u^{\nu}(x,0))-\ve^{p,0}(x))\id\,,
\end{eqnarray*}
which implies that $\dev\,(T^{\nu}_E(0))=\dev\,(T^0_E)$ and $b^{\nu}(0)=b^0$. From the assumption (\ref{eq:3.7}) 
we get that the last term on the right hand side of (\ref{eq:5.8}) is equal to zero.\\  
We know that the function $\|\nabla u\|_{L^2(\Omega;\R^9)}+\int\nolimits_{\Gamma_D}|u|dS$ is a norm on $H^1(\Omega;\R^3)$ equivalent to the standard norm (see \cite{21}). Applying this fact in (\ref{eq:5.8}) and the following well-known estimate [\cite{15}, p.36]\\
$$\|\nabla u\|^2_{L^2(\Omega)}\leq C^{\mathrm{curl}}_{\mathrm{div}}(\|\mathrm{div}\; u\|^2_{L^2(\Omega)}+
\|\mathrm{curl}\; u\|^2_{L^2(\Omega)}),$$
(the constant $C^{\mathrm{curl}}_{\mathrm{div}}$ does not depend on $u$ and 'curl' is the rotation operator) we get the following inequality 
\begin{eqnarray}
\label{eq:5.9}
&&\int\nolimits_{\Omega}\frac{1}{2\nu}\{|\dev\,(T_E^{\nu})(t)-b^{\nu}(t)|-\KK\}_{+}^2dx+\int\nolimits^t_0\int\nolimits_{\Omega}\D^{-1}T_{E,t}^{\nu}(\tau)T_{E,t}^{\nu}(\tau)dxd\tau\nn\\[1ex]
&+& 2\mu_c\int\nolimits^t_0\int\nolimits_{\Omega}|\mathrm{skew}(\nabla_x u_t^{\nu}(\tau))-A_t^{\nu}(\tau)|^2dxd\tau
+4l_c\int\nolimits^t_0\int\nolimits_{\Omega}|\nabla\mathrm{axl}(A_t^{\nu}(\tau))|^2dxd\tau\nn\\[1ex]
&\leq & \alpha\int\nolimits^t_0\Big(\|T_t^{\nu}(\tau)\|^2_{L^2(\Omega;\S)} + \|\mathrm{skew}(\nabla u^{\nu}_t(\tau))-A^{\nu}_t(\tau)\|^2_{L^2(\Omega;\R^9)} + \|\nabla\mathrm{axl}(A^{\nu}_t(\tau))\|^2_{L^2(\Omega;\R^9)}\Big)d\tau\nn\\[1ex]
&+&\alpha\int\nolimits^t_0\Big(\|\mathrm{div}\; u_t^{\nu}(\tau)\|^2_{L^2(\Omega;\R)}+ \|\mathrm{curl}\; u_t^{\nu}(\tau)\|^2_{L^2(\Omega;\R^3)}\Big)d\tau+\tilde{C}(T)\,.
\end{eqnarray}
From the observation $\mathrm{div}\;u^{\nu}_t=\mathrm{tr}\;\ve(u^{\nu}_t)-\mathrm{tr}\;\ve^{p,\nu}_t$ we have
\begin{eqnarray}
\label{eq:5.10}
\int\nolimits^t_0\|\mathrm{div}\; u_t^{\nu}(\tau)\|^2_{L^2(\Omega;\R)}d\tau\leq C\int\nolimits^t_0\|T_t^{\nu}(\tau)\|^2_{L^2(\Omega;\S)}d\tau\,,
\end{eqnarray}
where the constant $C$ does not depends on $\nu$. Now we choose $\alpha>0$ sufficiently small and we arrive at the inequality
\begin{eqnarray}
\label{eq:5.11}
&&\int\nolimits_{\Omega}\frac{1}{2\nu}\{|\dev\,(T_E^{\nu})(t)-b^{\nu}(t)|-\KK\}_{+}^2dx+\int\nolimits^t_0\int\nolimits_{\Omega}\D^{-1}T_{E,t}^{\nu}(\tau)T_{E,t}^{\nu}(\tau)dxd\tau\nn\\[1ex]
&+& 2\mu_c\int\nolimits^t_0\int\nolimits_{\Omega}|\mathrm{skew}(\nabla_x u_t^{\nu}(\tau))-A_t^{\nu}(\tau)|^2dxd\tau
+4l_c\int\nolimits^t_0\int\nolimits_{\Omega}|\nabla\mathrm{axl}(A_t^{\nu}(\tau))|^2dxd\tau\nn\\[1ex]
&\leq & \alpha\int\nolimits^t_0\|\mathrm{curl}\; u_t^{\nu}(\tau)\|^2_{L^2(\Omega;\R^3)}d\tau+\overline{C}(T,\alpha)\,.
\end{eqnarray}
Notice that
\begin{eqnarray}
\label{eq:5.12}
&&2\mu_c\int\nolimits^t_0\int\nolimits_{\Omega}|\mathrm{skew}(\nabla_x u_t^{\nu}(\tau))-A_t^{\nu}(\tau)|^2dxd\tau
+4l_c\int\nolimits^t_0\int\nolimits_{\Omega}|\nabla\mathrm{axl}(A_t^{\nu}(\tau))|^2dxd\tau\nn\\[1ex]
&=& 2\mu_c\int\nolimits^t_0\int\nolimits_{\Omega}|\mathrm{skew}(\nabla_x u_t^{\nu}(\tau))|^2dxd\tau+ 2\mu_c\int\nolimits^t_0\int\nolimits_{\Omega}|A_t^{\nu}(\tau)|^2dxd\tau\nn\\[1ex]
&-&4\mu_c\int\nolimits^t_0\int\nolimits_{\Omega}\mathrm{skew}(\nabla_x u_t^{\nu}(\tau))\cdot A_t^{\nu}(\tau)dxd\tau
+4l_c\int\nolimits^t_0\int\nolimits_{\Omega}|\nabla\mathrm{axl}(A_t^{\nu}(\tau))|^2dxd\tau\nn\\[1ex]
&\geq&\mu_c\int\nolimits^t_0\int\nolimits_{\Omega}|\mathrm{skew}(\nabla_x u_t^{\nu}(\tau))|^2dxd\tau - 
2\mu_c\int\nolimits^t_0\int\nolimits_{\Omega}|A_t^{\nu}(\tau)|^2dxd\tau\nn\\[2ex]
&+&4l_c\int\nolimits^t_0\int\nolimits_{\Omega}|\nabla\mathrm{axl}(A_t^{\nu}(\tau))|^2dxd\tau
\geq \qquad(\textrm{ Poincar\'e's inequality  })\qquad\geq\nn\\[2ex]
&\geq &\mu_c\int\nolimits^t_0\int\nolimits_{\Omega}|\mathrm{curl}\;u_t^{\nu}(\tau))|^2dxd\tau - 4\mu_c\,C_{\Omega}\int\nolimits^t_0\int\nolimits_{\Omega}|\nabla\mathrm{axl} A_t^{\nu}(\tau)|^2dxd\tau\nn\\[1ex]
&+&4l_c\int\nolimits^t_0\int\nolimits_{\Omega}|\nabla\mathrm{axl}(A_t^{\nu}(\tau))|^2dxd\tau\,,
\end{eqnarray}
where the constant $C_{\Omega}>0$ depends on the domain $\Omega$ only. Shifting the negative term on the right hand side of (\ref{eq:5.12}) into the left we may obtain the inequality
\begin{eqnarray}
\label{eq:5.13}
&&2\mu_c\int\nolimits^t_0\int\nolimits_{\Omega}|\mathrm{skew}(\nabla_x u_t^{\nu}(\tau))-A_t^{\nu}(\tau)|^2dxd\tau
+4l_c\int\nolimits^t_0\int\nolimits_{\Omega}|\nabla\mathrm{axl}(A_t^{\nu}(\tau))|^2dxd\tau\nn\\[1ex] 
&\geq & C(\Omega,\mu_c,l_c)\Big(\mu_c\int\nolimits^t_0\int\nolimits_{\Omega}|\mathrm{curl}\;u_t^{\nu}(\tau))|^2dxd\tau
+4l_c\int\nolimits^t_0\int\nolimits_{\Omega}|\nabla\mathrm{axl}(A_t^{\nu}(\tau))|^2dxd\tau\Big).\qquad\quad
\end{eqnarray}
If we apply inequality (\ref{eq:5.13}) to (\ref{eq:5.11}) then we can again choose $\alpha>0$ sufficiently small to end the proof.$\mbox{}$ \hfill $\Box$\\[2ex]
\begin{tw}
\label{tw:5.2}
Let us suppose that all hypotheses of Theorem \ref{tw:5.1} hold. Then the sequence $\{(\ve^{p,\nu}_t,b^{\nu}_t)\}_{\nu>0}$ is bounded in $L^2(0,T;L^2(\Omega;\SS)\times L^2(\Omega;\SS))$.
\end{tw}
{\bf\em Proof:}\hspace{2ex} The proof of the Theorem \ref{tw:5.1} yields that the sequence $\{u^{\nu}_t\}_{\nu>0}$ is bounded in $L^2(0,T; H^1(\Omega;\R^3))$. From the elastic constitutive equation we obtain
\begin{eqnarray}
\label{eq:5.14}
\|\ve^{p,\nu}_t(t)\|_{L^2(\Omega;\SS)}&\leq& \|\ve(u^{\nu}_t(t))\|_{L^2(\Omega;\S)}+\|\D^{-1}T^{\nu}_t(t)\|_{L^2(\Omega;\S)}\nn\\
&\leq& \|\ve(u^{\nu}_t(t))\|_{L^2(\Omega;\S)}+C\|T^{\nu}_t(t)\|_{L^2(\Omega;\S)},
\end{eqnarray}
where $C>0$ does not depend on $\nu>0$. Theorem \ref{tw:5.1} gives us the boundedness of the sequence $\{\ve^{p,\nu}_t\}_{\nu>0}$ in $L^2(0,T;L^2(\Omega;\SS))$. Next using equation $(\ref{eq:4.1})_6$ we have
\begin{eqnarray}
\label{eq:5.15}
\|b^{\nu}_t(t)\|^2_{L^2(\Omega;\SS)}&\leq & \tilde{C}\Big(c^2\,\|\ve^{p,\nu}_t(t)\|^2_{L^2(\Omega;\SS)}+d^2\,\|\,|\ve^{p,\nu}_t(t)|b^{\nu}(t)\|^2_{L^2(\Omega;\SS)}\Big)
\end{eqnarray}
and
\begin{eqnarray}
\label{eq:5.16}
\||\ve^{p,\nu}_t(t)|b^{\nu}(t)\|^2_{L^2(\Omega;\SS)}&=&\int\nolimits_{\Omega}|\ve^{p,\nu}_t(t)|^2|b^{\nu}(t)|^2dx\nn\\[1ex]
&\leq&\|\ve^{p,\nu}_t(t)\|^2_{L^2(\Omega;\SS)}\||b^{\nu}(t)|^2\|_{L^{\infty}(\Omega;\SS)}\,.
\end{eqnarray}
From Lemma \ref{lem:4.2} we know that $|b^{\nu}(x,t)|\leq \frac{c}{d}$ for almost all $(x,t)\in\Omega\times (0,T)$, then the term on the right hand side of (\ref{eq:5.16}) is bounded for almost all $t\in (0,T)$. This finishes the proof.$\mbox{}$ \hfill $\Box$\\[2ex]
{\bf\Large\em Proof of Theorem \ref{tw:3.3}}\\[1ex]
The motivation of the Definition \ref{de:3.2} yields that $L^2$-strong solutions of the system (\ref{eq:4.1}) with boundary and initial conditions (\ref{eq:4.2}) and (\ref{eq:4.3}) satisfy the energy inequality, which means that for all test functions $(\hat{T}_E,\hat{b})\in L^{2}(0,T;L^{2}(\Omega;\S)\times L^{2}(\Omega;\SS))$ such that 
$$(\dev\,(\hat{T}_E),\hat{b})\in \K^{\ast},\quad \mathrm{div}\; \hat{T}_E\in L^{2}(0,T;L^2(\Omega,\R^3))\,,\quad \hat{T}_E\cdot n=g_N\; \mathrm{on}\; \Gamma_N\times (0,T)$$
the inequality
\begin{eqnarray}
\label{eq:5.17}
&&\frac{1}{2}\int\nolimits_{\Omega}\D^{-1}T_{E}^{\nu}(x,t)T_E^{\nu}(x,t)dx+ \mu_c\int\nolimits_{\Omega}|\mathrm{skew}(\nabla_x u^{\nu}(x,t))-A^{\nu}(x,t)|^2dx\nn\\[1ex]
&+&2l_c\int\nolimits_{\Omega}|\nabla\mathrm{axl}(A^{\nu}(x,t))|^2dx+ \frac{1}{2c}\int\nolimits_{\Omega}|b^{\nu}(x,t)|^2dx
\leq\frac{1}{2}\int\nolimits_{\Omega}\D^{-1}T^0_{E}(x)T^0_E(x)dx\nn\\[1ex]
&+& \mu_c\int\nolimits_{\Omega}|\mathrm{skew}(\nabla_x u(x,0))-A(x,0)|^2dx+ \frac{1}{2c}\int\nolimits_{\Omega}|b(x,0)|^2dx+2l_c\int\nolimits_{\Omega}|\nabla\mathrm{axl}(A(x,0))|^2dx\nn\\[1ex]
 &+&\int\nolimits_0^t\int\nolimits_{\Omega}u_t^{\nu}(x,\tau)f(x,\tau)dxd\tau+ \int\nolimits_0^t\int\nolimits_{\Omega}u_t^{\nu}(x,\tau)\mathrm{div}\hat{T}_E(x,\tau))dxd\tau\\[1ex]
&+&\int\nolimits_0^t\int\nolimits_{\Gamma_D}g_{D,t}(x,\tau)(T^{\nu}(x,\tau)-\hat{T}_E(x,\tau))\cdot n(x)dSd\tau+ \int\nolimits_0^t\int\nolimits_{\Omega}\D^{-1}T_{E,t}^{\nu}(x,\tau)\hat{T}_E(x,\tau)dxd\tau\nn\\[1ex]
&+& \frac{1}{c}\int\nolimits_0^t\int\nolimits_{\Omega}b_t^{\nu}(x,\tau)\hat{b}(x,\tau)dxd\tau 
+4l_c\int\nolimits_0^t\int\nolimits_{\partial\Omega}\nabla\mathrm{axl}(A^{\nu}(x,\tau))\cdot n \;\mathrm{axl}(A_{D,t}(x,\tau))dSd\tau\nn
\end{eqnarray}
holds for all $t\in(0,T)$, where $T^0_E\in L^2(\Omega;\S)$ and $(u(0),A(0))\in H^{1}(\Omega;\R^3)\times H^2(\Omega;\mathfrak{so}(3))$ are unique solution of the problem (\ref{eq:3.8}). From Theorem \ref{tw:5.1} and \ref{tw:5.2} we conclude that for a subsequence (again denoted by $\nu$) we have 
\begin{eqnarray*}
u^{\nu}\rightharpoonup u\quad &\mathrm{in}& \quad L^{\infty}(0,T;H^1(\Omega;\R^3)),\\
u^{\nu}_t\rightharpoonup u_t\quad &\mathrm{in}& \quad L^{2}(0,T;H^1(\Omega;\R^3)),\\
A^{\nu}\rightharpoonup A\quad &\mathrm{in}& \quad L^{\infty}(0,T;H^2(\Omega;\so(3))),\\
T^{\nu}\rightharpoonup T\quad &\mathrm{in}& \quad L^{\infty}(0,T;L^2(\Omega;\R^9)),\\
T^{\nu}_t\rightharpoonup T_t\quad &\mathrm{in}& \quad L^{2}(0,T;L^2(\Omega;\R^9)),\\
b^{\nu}\rightharpoonup b\quad &\mathrm{in}& \quad L^{\infty}(0,T;L^2(\Omega;\SS)),\\
b^{\nu}_t\rightharpoonup b_t\quad &\mathrm{in}& \quad L^{2}(0,T;L^2(\Omega;\SS)).
\end{eqnarray*}
The right hand side of (\ref{eq:5.17}) has the limit for $\nu\rightarrow 0^+$, the left one is a convex functional thus, the sequential weak lower semicontinuity implies that
\begin{eqnarray}
\label{eq:5.18}
&&\frac{1}{2}\int\nolimits_{\Omega}\D^{-1}T_{E}(x,t)T_E(x,t)dx+ \mu_c\int\nolimits_{\Omega}|\mathrm{skew}(\nabla_x u(x,t))-A(x,t)|^2dx\nn\\[1ex]
&+&2l_c\int\nolimits_{\Omega}|\nabla\mathrm{axl}(A(x,t))|^2dx+ \frac{1}{2c}\int\nolimits_{\Omega}|b(x,t)|^2dx\nn\\[1ex]
&\leq&\liminf_{\nu\rightarrow 0^+}\Bigg\{\frac{1}{2}\int\nolimits_{\Omega}\D^{-1}T_{E}^{\nu}(x,t)T_E^{\nu}(x,t)dx+ \mu_c\int\nolimits_{\Omega}|\mathrm{skew}(\nabla_x u^{\nu}(x,t))-A^{\nu}(x,t)|^2dx\nn\\[1ex]
&+&2l_c\int\nolimits_{\Omega}|\nabla\mathrm{axl}(A^{\nu}(x,t))|^2dx+ \frac{1}{2c}\int\nolimits_{\Omega}|b^{\nu}(x,t)|^2dx\Bigg\}.
\end{eqnarray} 
Using (\ref{eq:5.18}) we complete the proof of the Theorem \ref{tw:3.3}.$\mbox{}$ \hfill $\Box$\\[6ex]
{\bf Acknowledgments: }This work has been supported by the European Union in
the framework of the European Social Fund through the Warsaw University of Technology Development Programme.

\footnotesize{
}
\end{document}